\documentclass[12]{article}

\usepackage{latexsym}
\usepackage{amsfonts}
\usepackage{graphicx}
\usepackage{amsmath}
\usepackage{amssymb}
\usepackage{float}
\usepackage{alltt}
\usepackage{eufrak}
\usepackage[latin1]{inputenc}

\makeatletter \@addtoreset{equation}{section}

\makeatletter

\def\be   {\begin{equation}}   \def\ee   {\end{equation}}
\def\ba   {\begin{array}}      \def\ea   {\end{array}}
\def\bea  {\begin{eqnarray}}   \def\eea  {\end{eqnarray}}
\def\bean {\begin{eqnarray*}}  \def\eean {\end{eqnarray*}}

\newtheorem{theorem} {Theorem}
\newtheorem{lemma}{Lemma}
\newtheorem{definition} {Definition}

\newtheorem{corollary} {Corollary}
\newtheorem{remark}{Remark}
\newtheorem{example} {Example}

\newtheorem{proposition} {Proposition}

\newcommand{\pre}{\mathrm{Re}}
\newcommand{\pim}{\mathrm{Im}}

\newcommand{\bi} {\ensuremath{{\bf i}}}
\newcommand{\bo} {\ensuremath{{\bf i_1}}}
\newcommand{\bos}{\ensuremath{{\bf i_1^{\text 2}}}}
\newcommand{\bt} {\ensuremath{{\bf i_2}}}
\newcommand{\bts}{\ensuremath{{\bf i_2^{\text 2}}}}

\newcommand {\bj}{\ensuremath{{\bf j}}}

\newcommand {\bjs}{\ensuremath{{\bf j^{\text 2}}}}
\newcommand{\eo} {\ensuremath{{\bf e_1}}}
\newcommand{\et} {\ensuremath{{\bf e_2}}}

\newcommand{\mC}{\ensuremath{\mathbb{C}}}

\newcommand{\mR}{\ensuremath{\mathbb{R}}}

\begin{document}

\vspace{4cm}
\begin{center} \LARGE{\textbf{Normal Families of Bicomplex Holomorphic Functions}}
\end{center}
\vspace{1cm}

\begin{center} \bf{K.S. Charak$^{1 }$,\quad D. Rochon$^{2}$,\quad N. Sharma$^{3}$ }
\end{center}

\bigskip

\begin{center}
{$^{1}$ Department of Mathematics, University of Jammu,
Jammu-180 006, INDIA.\\
E-mail: kscharak7@rediffmail.com }
\end{center}

\medskip
\begin{center} {$^{2}$ D\'epartement de math\'ematiques et
d'informatique, Universit\'e du Qu\'ebec \`a Trois-Rivi\`eres, C.P. 500 Trois-Rivi\`eres, Qu\'ebec, Canada G9A 5H7.  \\
E-mail: Dominic.Rochon@UQTR.CA, Web: www.3dfractals.com}
\end{center}

\medskip

\begin{center}{$^{3}$ Department of Mathematics, University of Jammu,
Jammu-180 006, INDIA.\\
E-mail: narinder25sharma@sify.com }
\end{center}

\bigskip
\begin{abstract}
In this article, we introduce the concept of normal families of
bicomplex holomorphic functions to obtain a bicomplex Montel
theorem. Moreover, we give a general definition of Fatou and Julia
sets for bicomplex polynomials and we obtain a characterization of
bicomplex Fatou and Julia sets in terms of Fatou set, Julia set
and filled-in Julia set of one complex variable. Some 3D visual
examples of bicomplex Julia sets are also given for the specific
slice $\bj=0$.
\end{abstract}

\vspace{1cm} \noindent \textbf{Keywords: }Bicomplex Numbers,
Complex Clifford Algebras, Normal Families, Montel Theorem, Bicomplex Dynamics, Julia Set, Fatou Set.\\

\vspace{4cm}

\normalsize
\newpage

\section{Introduction}

\bigskip

A family $\boldsymbol F$ of holomorphic functions defined on a
domain $D \subseteq \mC$ is said to be normal in $D$ if every
sequence in $\boldsymbol F$ has a subsequence converging uniformly
on compact subsets of $D$ to a function $f$. The limit function $f$
is holomorphic on $D$ (by Weierstrass Theorem) or the constant
infinity. Various authors while studying the normality of a family
of holomorphic functions take the limit function $f\neq \infty$
but for studying the normal families from the complex dynamics point
of view, one needs to include the case where the limit function
$f\equiv \infty$. The former approach we shall call as restrictive
approach while the later will be called the general approach
towards normal families. The concept of normal families was
introduced by P. Montel in 1907 \cite{19}. For a comprehensive
account of normal families of meromorphic functions on domains in
$\mC$ one can refer to Joel Schiff's text \cite{27}, C.T. Chuang's
text \cite{20} and Zalcman's survey article \cite{21}. With the
renewed interest in normal families of meromorphic functions,
arising largely from the important role they play in Complex
Dynamics, it seems sensible to talk about normal families of
holomorphic functions on different domains of different spaces
thereby enabling one to study the dynamics of such functions. In
this article we have considered the families of bicomplex
holomorphic functions on bicomplex domains. Since this article
lays the foundations of the subject {\bf Normal Families of
Bicomplex Holomorphic Functions} for future investigations in
various possible directions, it is necessary to adopt a dual
approach towards the study of normality of families of bicomplex
holomorphic functions on bicomplex domains. The first approach
is a restrictive approach which gives rise to more interesting
results when the normal families are studied in their own right.
For example the converse of the Montel Theorem holds under this
approach. The second approach is the general approach in which
though the converse of the Montel Theorem fails to hold but is
essentially required when the normal families are studied from the
bicomplex dynamics point of view. During our discussions, we shall
come across the situations where the differences lead to
interesting conclusions. Besides a complete discussion on the Montel
Theorem in various situations, we have defined Fatou, Julia and
filled-in Julia sets of bicomplex polynomials and their
representations in terms of their complex counterparts are
obtained for the specific case of non-degenerate bicomplex
polynomials of degree $d\geq 2$. Also, some 3D visual examples of
bicomplex Julia sets are given for the specific slice $\bj=0$.

\section{Preliminaries}

\subsection{Bicomplex Numbers}
Bicomplex numbers are defined as
\be \mathbb{T}:=\{z_1+z_2\bt\ |\ z_1, z_2 \in \mathbb{C}(\bo) \}
\label{enstetra}
\ee where the imaginary units $\bo, \bt$ and $\bj$
are governed by the rules: $\bos=\bts=-1$, $\bjs=1$ and \be
\ba{rclrcl}
\bo\bt &=& \bt\bo &=& \bj,  \\
\bo\bj &=& \bj\bo &=& -\bt, \\
\bt\bj &=& \bj\bt &=& -\bo. \\
\ea \ee Note that we define $\mC(\bi_k):=\{x+y\bi_k\ |\ \bi_k^2= -1$
and $x,y\in \mR \}$ for $k=1,2$. Hence, it is easy to see that the
multiplication of two bicomplex numbers is commutative. In fact, the bicomplex numbers
$$\mathbb{T}\cong {\rm Cl}_{\Bbb{C}}(1,0) \cong {\rm Cl}_{\Bbb{C}}(0,1)$$
are \textbf{unique} among the \textit{Complex Clifford algebras} (see \cite{32})
in that they are commutative but not division algebra.
It is also convenient to write the set of bicomplex numbers as
\be
\mathbb{T}:=\{w_0+w_1\bo+w_2\bt+w_3\bj\ |\ w_0,
 w_1,w_2,w_3 \in \mathbb{R}\}.
\ee

In particular, in the equation (\ref{enstetra}), if we put $z_1=x$ and
$z_2=y\bo$ with $x,y \in \mathbb{R}$, then we obtain the following
subalgebra of hyperbolic numbers, also called duplex numbers (see, e.g. \cite {18}, \cite {24}):
$$\mathbb{D}:=\{x+y\bj\ |\ \bj^2=1,\ x,y\in \mathbb{R}\}\cong {\rm Cl}_{\Bbb{R}}(0,1).$$

Complex conjugation plays an important role both for algebraic and
geometric properties of $\mathbb{C}$. For bicomplex numbers, there are three possible
conjugations. Let $w\in \mathbb{T}$ and $z_1,z_2 \in
\mathbb{C}(\mathbf{i_1})$ such that $w=z_1+z_2\mathbf{i_2}$. Then we
define the three conjugations as:

\begin{subequations}
\label{eq:dag}
\begin{align}
w^{\dag_{1}}&=(z_1+z_2\bt)^{\dag_{1}}:=\overline z_1+\overline z_2
\bt,
\\
w^{\dag_{2}}&=(z_1+z_2\bt)^{\dag_{2}}:=z_1-z_2 \bt,
\\
w^{\dag_{3}}&=(z_1+z_2\bt)^{\dag_{3}}:=\overline z_1-\overline z_2
\bt,
\end{align}
\end{subequations}
where $\overline z_k$ is the standard complex conjugate of complex
number $z_k \in \mathbb{C}(\mathbf{i_1})$. If we say that the
bicomplex number $w=z_1+z_2\bt=w_0+w_1\bo+w_2\bt+w_3\bj$ has the
``signature'' $(++++)$, then the conjugations of type 1,2 or 3 of
$w$ have, respectively, the signatures $(+-+-)$, $(++--)$ and
$(+--+)$. We can verify easily that the composition of the
conjugates gives the four-dimensional Klein group:
\begin{center}
\be
\begin{tabular}{|c||c|c|c|c|}
\hline
$\circ$ & $\dag_{0}$ & $\dag_{1}$  & $\dag_{2}$  & $\dag_{3}$  \\
\hline
\hline
$\dag_{0}$    & $\dag_{0}$ & $\dag_{1}$  & $\dag_{2}$  & $\dag_{3}$  \\
\hline
$\dag_{1}$    & $\dag_{1}$ & $\dag_{0}$ & $\dag_{3}$  & $\dag_{2}$ \\
\hline
$\dag_{2}$    & $\dag_{2}$ & $\dag_{3}$  & $\dag_{0}$ & $\dag_{1}$ \\
\hline
$\dag_{3}$    & $\dag_{3}$ & $\dag_{2}$ & $\dag_{1}$ & $\dag_{0}$ \\
\hline
\end{tabular}
\label{eq:groupedag} \ee
\end{center}
where $w^{\dag_{0}}:=w\mbox{ } \forall w\in \mathbb{T}$.

The three kinds of conjugation all have some of the standard properties of
conjugations, such as:
\begin{eqnarray}
(s+ t)^{\dag_{k}}&=&s^{\dag_{k}}+ t^{\dag_{k}},\\
\left(s^{\dag_{k}}\right)^{\dag_{k}}&=&s, \\
\left(s\cdot t\right)^{\dag_{k}}&=&s^{\dag_{k}}\cdot t^{\dag_{k}},
\end{eqnarray}
for $s,t \in \mathbb{T}$ and $k=0,1,2,3$.\\

We know that the product of a standard complex number with its
conjugate gives the square of the Euclidean metric in
$\mathbb{R}^2$. The analogs of this, for bicomplex numbers, are
the following. Let $z_1,z_2 \in \mathbb{C}(\bo)$ and
$w=z_1+z_2\bt\in \mathbb{T}$, then we have that :
\begin{subequations}
\begin{align}
|w|^{2}_{\bo}&:=w\cdot w^{\dag_{2}}=z^{2}_{1}+z^{2}_{2} \in
\mathbb{C}(\bo),
\\*[2ex] |w|^{2}_{\bt}&:=w\cdot
w^{\dag_{1}}=\left(|z_1|^2-|z_2|^2\right)+2\pre(z_1\overline
z_2)\bt \in \mathbb{C}(\bt), \\*[2ex] |w|^{2}_{\bj}&:=w\cdot
w^{\dag_{3}}=\left(|z_1|^2+|z_2|^2\right)-2\pim(z_1\overline
z_2)\bj \in \mathbb{D},
\end{align}
\end{subequations}
where the subscript of the square modulus refers to the subalgebra
$\mathbb{C}(\bo), \mathbb{C}(\bt)$ or $\mathbb{D}$ of $\mathbb{T}$
in which $w$ is projected. Note that for $z_1,z_2 \in \mathbb{C}(\bo)$ and $w=z_1+z_2\bt\in
\mathbb{T}$, we can define the usual (Euclidean in $\mR^4$) norm
of $w$ as $\|w\|=\sqrt{|z_1|^2+|z_2|^2}=\sqrt{\pre(|w|^{2}_{\bj})}$.

It is easy to verify that $w\cdot \displaystyle
\frac{w^{\dag_{2}}}{|w|^{2}_{\bo}}=1$. Hence, the inverse of $w$
is given by \be w^{-1}= \displaystyle
\frac{w^{\dag_{2}}}{|w|^{2}_{\bo}}. \ee From this, we find that
the set $\mathcal{NC}$ of zero divisors of $\mathbb{T}$, called
the {\em null-cone}, is given by $\{z_1+z_2\bt\ |\
z_{1}^{2}+z_{2}^{2}=0\}$, which can be rewritten as \be
\mathcal{NC}=\{z(\bo\pm\bt)|\ z\in \mathbb{C}(\bo)\}. \ee
\smallskip

\subsection{Bicomplex Holomorphic Functions}

It is also possible to define differentiability of a function at a
point of $\mathbb{T}$:

\begin{definition}
Let $U$ be an open set of $\mathbb{T}$ and $w_0\in U$. Then,
$f:U\subseteq\mathbb{T}\longrightarrow\mathbb{T}$ is said to be
$\mathbb{T}$-differentiable at $w_{0}$ with derivative equal to
$f^\prime(w_0)\in\mathbb{T}$ if
$$\lim_{\stackrel{\scriptstyle w \rightarrow w_{0}}
{\scriptscriptstyle (w-w_{0}\mbox{
}inv.)}}\frac{f(w)-f(w_{0})}{w-w_{0}} =f^\prime(w_0).$$
\textbf{Note:} The subscript of the limit is there to recall that division is possible only if $w-w_{0}$
is invertible.
\end{definition}

We also say that the function $f$ is bicomplex holomorphic ($\mathbb{T}$-holomorphic) on
an open set U if and only if $f$ is $\mathbb{T}$-differentiable at
each point of U.

\smallskip
Using $w=z_1+z_2\bt$, the bicomplex number $w$ can be seen as an element $(z_1,z_2)$ of
$\mathbb{C}^{2}$, so a function
$f(z_1+z_2\bold{i_2})=f_1(z_1,z_2)+f_2(z_1,z_2)\bold{i_2}$ of
$\mathbb{T}$ can be seen as a mapping
$f(z_1,z_2)=(f_1(z_1,z_2),f_2(z_1,z_2))$ of $\mathbb{C}^{2}$. Here
we have a characterization of such mappings:

\begin{theorem}
Let $U$ be an open set and
$f:U\subseteq\mathbb{T}\longrightarrow\mathbb{T}$ such that $f\in
{C}^{1}(U)$, and let
$f(z_1+z_2\bold{i_2})=f_1(z_1,z_2)+f_2(z_1,z_2)\bold{i_2}$. Then
$f$ is $\mathbb{T}$-holomorphic on $U$ if and only if
$$\mbox{$f_1$ and $f_2$ are holomorphic in $z_1$ and $z_2$},$$
and
$$\frac{\partial{f_1}}{\partial{z_1}}
=\frac{\partial{f_2}}{\partial{z_2}}\mbox{ and
}\frac{\partial{f_2}}{\partial{z_1}}
=-\frac{\partial{f_1}}{\partial{z_2}}\mbox{ on U}.$$
\\
Moreover, $f^\prime=\frac{\partial{f_1}}{\partial{z_1}}
+\frac{\partial{f_2}}{\partial{z_1}}\bold{i_2}$ and $f^\prime(w)$ is
invertible if and only if  $det\mathcal{J}_{f}(w)\neq 0$ where $\mathcal{J}_{f}$ is the Jacobian matrix of $f$.
\label{theobasic}
\end{theorem}

\smallskip
This theorem can be obtained from the results in \cite{14} and
\cite{17}. Moreover, by the Hartogs theorem \cite{23}, it is
possible to show that ``$f\in C^{1}(U)$" can be dropped from the
hypotheses. Hence, it is natural to define the corresponding class
of mappings for $\mathbb{C}^{2}$:

\begin{definition}
The class of $\mathbb{T}$-holomorphic mappings on a open set
$U\subseteq\mathbb{C}^{2}$ is defined as follows:
$$
\mbox{$TH(U):=$}\{f\mbox{:$U$}\subseteq\mathbb{C}^{2}\longrightarrow\mathbb{C}^{2}
|f\in\mbox{H($U$) and }\frac{\partial{f}_1}{\partial{z}_{1}}=
\frac{\partial{f}_2}{\partial{z}_{2}}\mbox{,
}\frac{\partial{f}_2}{\partial{z}_{1}}=-
\frac{\partial{f}_1}{\partial{z}_{2}}\mbox{ on $U$}\}.
$$

It is the subclass of holomorphic mappings of $\mathbb{C}^{2}$
satisfying the complexified Cauchy-Riemann equations.
\end{definition}

We remark that $f\in TH(U)$ in terms of $\mathbb{C}^{2}$ if and only if
$f$ is $\mathbb{T}$-differentiable on $U$.
It is also important to know that
every bicomplex number $z_1+z_2\bold{i_2}$ has the following
unique idempotent representation:
\begin{equation}
z_1+z_2\bt=(z_1-z_2\bo)\eo+(z_1+z_2\bo)\et.
\label{idempotent}
\end{equation}
where $\bold{e_1}=\frac{1+\bold{j}}{2}$ and $\bold{e_2}=\frac{1-\bold{j}}{2}$.
This representation is very useful because addition, multiplication and division can be done term-by-term.
It is also easy to verify the following characterization of the non-invertible elements.

\begin{proposition}
An element $w=z_1+z_2\bt$ is in the null-cone if and only if $z_1-z_2\bold{i_1}=0$
or $z_1+z_2\bold{i_1}=0$.
\label{Null}
\end{proposition}

The notion of the holomorphicity can also be seen with this kind of
notation. For this we need to define the projections
$\mathcal{P}_1,\mathcal{P}_2:\mathbb{T}\longrightarrow\mathbb{C}(\bold{i_1})$ as
$\mathcal{P}_1(z_1+z_2\bold{i_2})=z_1-z_2\bold{i_1}$ and
$\mathcal{P}_2(z_1+z_2\bold{i_2})=z_1+z_2\bold{i_1}$.
\begin{definition}
We say that $X\subseteq\mathbb{T}$ is a $\mathbb{T}$-cartesian set
determined by $X_1$ and $X_2$ if $X=X_{1}\times_e
X_{2}:=\{z_1+z_2\bold{i_2}\in\mathbb{T}:z_1+z_2\bold{i_2}=w_1\bold{e_1}+w_2\bold{e_2},
(w_1,w_2)\in X_1\times X_2\}$.
\end{definition}

In \cite{14} it is shown that if $X_1$ and $X_2$ are domains (open and connected) of
$\mathbb{C}(\bold{i_1})$ then $X_1\times_e X_2$ is also a domain
of $\mathbb{T}$. Then, a way to construct some ``discus" (of
center 0) in $\mathbb{T}$ is to take the $\mathbb{T}$-cartesian
product of two discs (of center 0) in $\mathbb{C}(\bold{i_1})$.
Hence, we define the ``discus" with center $a=a_1+a_2\bt$ of
radius $r_1$ and $r_2$ of $\mathbb{T}$ as follows \cite{14}:
$D(a;r_1,r_2):=B^{1}(a_1-a_2\bo,r_1)\times_{e}B^{1}(a_1+a_2\bo,r_2)=
\{z_1+z_2\bold{i_2}:z_1+z_2\bold{i_2}=w_1\bold{e_1}+w_2\bold{e_2},|w_1-(a_1-a_2\bo)|<r_1,|w_2-(a_1+a_2\bo)|<r_2\}$
where $B^{1}(z,r)$ is an open ball with center
$z\in\mathbb{C}(\bold{i_1})$ and radius $r>0$. In
the particular case where $r=r_1=r_2$, $D(a;r,r)$ will be
called the $\mathbb{T}$-disc with center $a$ and radius
$r$. In particular, we define $\overline{D}(a;r_1,r_2):=\overline{B^{1}(a_1-a_2\bo,r_1)}\times_{e}\overline{B^{1}(a_1+a_2\bo,r_2)}\subset \overline{D(a;r_1,r_2)}$.
We remark that $D(0;r,r)$ is, in fact, the Lie Ball
(see \cite{25}) of radius $r$ in $\mathbb{T}$.

\smallskip\smallskip\smallskip\smallskip
\noindent Now, it is possible to state the following
striking theorems (see \cite{14}):

\begin{theorem}
Let $X_1$ and $X_2$ be open sets in $\mathbb{C}(\bold{i_1})$. If $f_{e1}:X_1\longrightarrow \mathbb{C}(\bold{i_1})$ and
$f_{e2}:X_2\longrightarrow \mathbb{C}(\bold{i_1})$ are holomorphic
functions of $\mathbb{C}(\bold{i_1})$ on $X_1$ and
$X_2$ respectively, then the function $f:X_1\times_e
X_2\longrightarrow \mathbb{T}$ defined as
$$f(z_1+z_2\bold{i_2})=f_{e1}(z_1-z_2\bold{i_1})\bold{e_1}+f_{e2}(z_1+z_2\bold{i_1})\bold{e_2} \mbox{ }\forall\mbox{ }z_1+z_2\bold{i_2}\in X_1\times_e X_2$$
is $\mathbb{T}$-holomorphic on the open set $X_1\times_e X_2$ and
$$f^\prime(z_1+z_2\bold{i_2})=f^\prime_{e1}(z_1-z_2\bold{i_1})\bold{e_1}+f^\prime_{e2}(z_1+z_2\bold{i_1})\bold{e_2}$$
$\mbox{ }\forall\mbox{ }z_1+z_2\bold{i_2}\in X_1\times_e X_2.$
\label{theo5}
\end{theorem}

\begin{theorem}
Let $X$ be an open set in $\mathbb{T}$, and let
$f:X\longrightarrow\mathbb{T}$ be a $\mathbb{T}$-holomorphic
function on X. Then there exist holomorphic functions
$f_{e1}:X_1\longrightarrow\mathbb{C}(\bold{i_1})$ and
$f_{e2}:X_2\longrightarrow\mathbb{C}(\bold{i_1})$ with
$X_1=\mathcal{P}_1(X)$ and $X_2=\mathcal{P}_2(X)$, such that:
$$f(z_1+z_2\bold{i_2})=f_{e1}(z_1-z_2\bold{i_1})\bold{e_1}+f_{e2}(z_1+z_2\bold{i_1})\bold{e_2}$$
$\forall z_1+z_2\bold{i_2}\in X.$
\label{theo4}
\end{theorem}

\section{Bicomplex Montel Theorem}

\noindent  Since the concepts like {\bf uniform boundedness}, {\bf local uniform boundedness}, and {\bf uniform convergence on compact sets} are defined for functions on any metric space and do not depend on bicomplex holomorphic functions, we assume these concepts in our discussion and refer the reader to any standard text on Analysis(e.g. see \cite {33} and \cite {22}). We start our discussion with the following definition of normality.

\begin{definition}  A family $\boldsymbol F$ of bicomplex holomorphic functions defined on a domain $D \subseteq \mathbb{T}$ is said to be
\textbf{normal} in $D$ if every sequence in  $\boldsymbol F$
contains a subsequence which converges locally uniformly on $D$.
$\boldsymbol F$ is said to be \textbf{normal at a point} $z\in D$
if it is normal in some neighborhood of $z$ in $D$.
\label{basicnormal}
\end{definition}

\noindent Let $f:D\longrightarrow\mathbb{T}$ be a $\mathbb{T}$-holomorphic
function on $D$. Then by Theorem \ref{theo4}, there exist
holomorphic functions
$f_{e1}:\mathcal{P}_1(D)\longrightarrow\mathbb{C}(\bold{i_1})$ and
$f_{e2}:\mathcal{P}_2(D)\longrightarrow\mathbb{C}(\bold{i_1})$
such that
$$f(z_1+z_2\bold{i_2})=f_{e1}(z_1-z_2\bold{i_1})\bold{e_1}+f_{e2}(z_1+z_2\bold{i_1})\bold{e_2} \mbox{ }\forall
\mbox{ }z_1+z_2\bold{i_2}\in D.$$ We define the \textbf{norm}
of $f$ on $D$ as
$$\|f\|=\|f(z)\|=\boldsymbol\{\frac{|f_{e1}(z_1-z_2\bold{i_1})|^{2}+|f_{e2}(z_1+z_2\bold{i_1})|^{2}}{2}
\boldsymbol\}^{\frac{1}{2}},\\ z=z_1+z_2\bold{i_2} \in D.$$ One
can easily see that
\begin{itemize}
  \item $ \|f\| \geq 0, \|f\|=0  \text{ iff }\ f\equiv 0 \text{ on }\ D$;
  \item  $\|af\|=|a|\|f\|, \ a\in \mathbb{C}(\bold{i_1})$;
  \item  $\|f+g\|\leq \|f\|+\|g\|$;
  \item $\|fg\|\leq \sqrt{2}\|f\| \|g\|$.
\end{itemize}

\noindent Thus, the linear space of bicomplex holomorphic functions
on a domain $D \subseteq \mathbb{T}$ is a normed space under the
above norm.

\medskip
We start with a \textbf{uniformly bounded} family $\boldsymbol F$
of bicomplex holomorphic functions. In this case, we can verify
directly the following result.

\begin{theorem}  A family $\boldsymbol F$ of bicomplex holomorphic functions defined on a bicomplex cartesian domain $D$ is uniformly bounded on $D$
if and only if $\boldsymbol F_{ei}=\mathcal{P}_{i}(\boldsymbol F)$
is uniformly bounded on $\mathcal{P}_{i}(D), \ i=1,2. $
\label{unibounded}
\end{theorem}

If we consider now a \textbf{locally uniformly bounded} family
$\boldsymbol F$ of bicomplex holomorphic functions, we can prove a
similar result since a set $K=\mathcal{P}_1(K) \times_e
\mathcal{P}_2(K)$ is compact if and only if $\mathcal{P}_i(K)$ is
compact for $i=1,2$.

\begin{theorem} A family $\boldsymbol F$ of bicomplex holomorphic functions defined on a bicomplex cartesian domain $D$
is locally uniformly bounded on $D$ if and only if $\boldsymbol F_{ei}=\mathcal{P}_{i}(\boldsymbol F)$ is locally uniformly bounded on $\mathcal{P}_{i}(D), \ i=1,2.$
\label{thmP}
\end{theorem}
\emph{Proof}  Let $\boldsymbol F$ be locally uniformly bounded on
$D$. Then for every compact set $K \subset D$ there is a constant
$M(K)$ such that
$$\|f(z)\|\leq M, \ \forall f\in \boldsymbol F, \\ z=z_1+z_2\bold{i_2} \in K.$$
Thus,
$$\boldsymbol\{\frac{|f_{e1}(z_1-z_2\bold{i_1})|^{2}+|f_{e2}(z_1+z_2\bold{i_1})|^{2}}{2}
\boldsymbol\}^{\frac{1}{2}}\leq M,\ \forall f_{ei}\in \boldsymbol
F_{ei}, \ i=1,2 $$
$$  \forall z_1-z_2\bold{i_1}\in \mathcal{P}_{1}(K),\ z_1+z_2\bold{i_1}\in \mathcal{P}_2(K).$$
Therefore,
\begin{equation}
|f_{e1}(z_1-z_2\bold{i_1})|\leq \sqrt{2}M,\ \forall f_{e1}\in
\boldsymbol F_{e1}, \ \forall z_1-z_2\bold{i_1}\in
\mathcal{P}_{1}(K) \label{eq04}
\end{equation}
and
\begin{equation}
|f_{e2}(z_1+z_2\bold{i_1})|\leq \sqrt{2}M, \ \forall f_{e2}\in
\boldsymbol F_{e2} \ \forall z_1+z_2\bold{i_1}\in
\mathcal{P}_{2}(K). \label{eq05}
\end{equation}

Now, let $K_1$ be a compact subset of $\mathcal{P}_{1}(D)$. Then
there is (always) a compact subset $K_2$ of $\mathcal{P}_{2}(D)$
(even singleton will do) such that $K_1\times_{e} K_2 =
K^{\prime}$ say, is a compact subset of $D$ with
$\mathcal{P}_{i}(K^{\prime})=K_{i}, \ i=1,2.$ Thus (\ref{eq04})
holds for any compact subset of $\mathcal{P}_{1}(D)$, and
similarly for (\ref{eq05}).

Conversely, suppose $\boldsymbol F_{ei}$ is locally uniformly
bounded  on $\mathcal{P}_i(D),\ i=1,2.$ Let $K$ be any compact
subset of $D$. Then by continuity of $\mathcal{P}_i$,
$K_i=\mathcal{P}_i(K)$ is compact subset of $\mathcal{P}_i(D),\
i=1,2$ and hence there are constants $M_1(K_1)$ and $M_2(K_2)$
such that
$$ |f_{e1}(z_1-z_2\bold{i_1})|\leq M_1,\ \forall f_{e1}\in \boldsymbol F_{e1}, \ \forall z_1-z_2\bold{i_1}\in K_{1}$$
and
$$|f_{e2}(z_1+z_2\bold{i_1})|\leq M_2, \ \forall f_{e2}\in \boldsymbol F_{e2} \ \forall z_1+z_2\bold{i_1}\in K_2.$$
Therefore,
\begin{equation}
\|f(z)\|\leq \boldsymbol
\{\frac{{M_{1}}^{2}+{M_{2}}^{2}}{2}\boldsymbol \}^{\frac{1}{2}}, \
\forall f\in \boldsymbol F, \ z=z_1+z_2\bold{i_2} \in
\mathcal{P}_1(K) \times_{e}  \mathcal{P}_2(K). \label{eq06}
\end{equation}
Since $K \subseteq \mathcal{P}_1(K) \times_{e}  \mathcal{P}_2(K)$,
(\ref{eq06}) holds for $K$ also and this completes the proof.
$\Box$

What happens if $D$ is not a bicomplex cartesian product? In the
case of uniformly bounded family of bicomplex holomorphic
functions (Theorem \ref{unibounded}), it is easy to verify that the
result is true for any domain. In the case of locally uniformly bounded family of bicomplex holomorphic
functions, we need to recall the following results from the bicomplex function
theory.

\begin{remark}
A domain $D\subseteq \mathbb{T}$ is a domain of holomorphism for
functions of a bicomplex variable if and only if $D$ is a
$\mathbb{T}$-cartesian set (\cite {14}, Theorem 15.11), and if $D$
is not a domain of holomorphism then $D\subsetneq
\mathcal{P}_1(D)\times_{e}\mathcal{P}_2(D)$, and there exists a
holomorphic function which is a bicomplex holomorphic continuation
of the given function from $D$ to
$\mathcal{P}_1(D)\times_{e}\mathcal{P}_2(D)$ (\cite {14}, Corollary
15.4). \label{rek}
\end{remark}

\begin{theorem}
A family $\boldsymbol F$ of bicomplex holomorphic functions defined on a arbitrary bicomplex domain $D$
is locally uniformly bounded on $D$ if and only if $\boldsymbol F_{ei}=\mathcal{P}_{i}(\boldsymbol F)$ is locally uniformly bounded on $\mathcal{P}_{i}(D), \ i=1,2.$
\label{lem01}
\end{theorem}
\emph{Proof}
If $\boldsymbol F_{ei}=\mathcal{P}_{i}(\boldsymbol F)$ is locally uniformly bounded on $\mathcal{P}_{i}(D)$ for $i=1,2$, from Remark \ref{rek},
we can extend $D$ to $\mathcal{P}_1(D)\times_{e}\mathcal{P}_2(D)$ and apply Theorem \ref{thmP} to obtain that $\boldsymbol F$
is locally uniformly bounded on $\mathcal{P}_1(D)\times_{e}\mathcal{P}_2(D)$.
For the other side, we need to recall that a family $\boldsymbol F$ is locally uniformly bounded on $D$ if and only if
the family $\boldsymbol F$ is \textbf{locally bounded} on $D$ i.e. for each $w_0\in D$ there is a positive number $M=M(w_0)$
and a neighbourhood $D(w_0;r,r)\subset D$ such that $||f(w)||\leq M$ for all $w\in D(w_0;r,r)$ and all $f\in\boldsymbol F$ (see \cite{27}).
Since $D(w_0;r,r)\subset D$ is a bicomplex cartesian product of two discs in the plane, it is easy to verify that
the family $\boldsymbol F_{ei}$ is bounded by $\sqrt{2}M(w_0)$ on $D(\mathcal{P}_i(w_0),r)\subset \mathcal{P}_i(D)$ for $i=1,2$.
As $w_0$ was arbitrary, $\boldsymbol F_{ei}=\mathcal{P}_{i}(\boldsymbol F)$ is locally bounded on $\mathcal{P}_{i}(D), \ i=1,2$. $\Box$

\medskip
\noindent We are now ready to prove the bicomplex version of the Montel theorem.

\begin{lemma}  Let $\boldsymbol F$ be a family of bicomplex holomorphic functions defined on a bicomplex domain $D$.
If $F_{ei}=\mathcal{P}_{i}(\boldsymbol F)$ is normal on $\mathcal{P}_{i}(D)$ for $i=1,2$
then $\boldsymbol F$ is normal on $D$.
\label{normal}
\end{lemma}
\emph{Proof}
Suppose that $\boldsymbol F_{ei}=\mathcal{P}_{i}(\boldsymbol F)$ is normal on $\mathcal{P}_{i}(D)=D_i, \ i=1,2.$
We want to show that $\boldsymbol F$ is normal in $D$. Let $\{F_n\}$ be any sequence in $\boldsymbol F$ and
$K$ be any compact subset of $D$. Then $\{\mathcal{P}_1(F_n)\}=\{(f_n)_1\}$ is a sequence in
$\boldsymbol F_{e1}=\mathcal{P}_{1}(\boldsymbol F).$ Since $\boldsymbol F_{e1}=\mathcal{P}_{1}(\boldsymbol F)$ is normal in
$\mathcal{P}_1(D)$ then $\{(f_n)_1\}$ has a subsequence $\{(f_{n_k})_1 \}$ which converges
uniformly on $\mathcal{P}_1(K)$ to a $\mathbb{C}(\bo)$-function. Now, consider $\{F_{n_k}\}$ in $\boldsymbol F$.
Then $\{\mathcal{P}_2(F_{n_k})\}=\{(f_{n_k})_2\}$ is a sequence in
$\boldsymbol F_{e2}=\mathcal{P}_{2}(\boldsymbol F)$. Since $\boldsymbol F_{e2}=\mathcal{P}_{2}(\boldsymbol F)$ is normal in
$\mathcal{P}_1(D)$ then $\{(f_{n_k})_2\}$ has a subsequence $\{(f_{n_{k_l}})_2 \}$ which converges
uniformly on $\mathcal{P}_2(K)$ to a $\mathbb{C}(\bo)$-function. This implies that $\{(f_{n_{k_l}})_1\eo+(f_{n_{k_l}})_2\et\}$
is a subsequence of $\{\boldsymbol F_n\}$ which converges uniformly on $\mathcal{P}_1(K)\times_{e}\mathcal{P}_2(K)\supseteq K$ to a
bicomplex function showing that $\boldsymbol F$ is normal in $D$. $\Box$

\begin{theorem}\textbf{(Montel)}
Every locally uniformly bounded family of bicomplex holomorphic
functions defined on a bicomplex domain is a normal family.
\label{Montel}
\end{theorem}
\emph{Proof}
Let $\boldsymbol F$ be a locally uniformly bounded family of bicomplex holomorphic functions defined on a domain $D\subseteq \mathbb{T}$.
Using Theorem \ref{lem01}, we have that $\boldsymbol F_{ei}=\mathcal{P}_{i}(\boldsymbol F)$ is locally uniformly bounded on $\mathcal{P}_{i}(D), \ i=1,2.$
Hence, from the classical Montel Theorem, $\boldsymbol F_{ei}=\mathcal{P}_{i}(\boldsymbol F)$ is normal on $\mathcal{P}_{i}(D)$ for  $i=1,2$
and by Lemma \ref{normal} we obtain that $\boldsymbol F$ is normal on $D$.$\Box$

\medskip
\textbf{Note:} The converse of the Bicomplex Montel Theorem is also
true. Indeed, suppose that $\boldsymbol F$ is normal and not
locally uniformly bounded in $D$. Then in some closed discus
$\overline{D}(a;r_1,r_2)$ in the domain $D$, for each $n\in\mathbb
N$ there is a function $f_n \in \boldsymbol F$ and a point $w_n
\in \overline{D}(a;r_1,r_2)$ such that $\|f_{n}(w_{n}) \|>n.$
Since $\boldsymbol F$ is normal, there is a subsequence
$\{f_{n_k}\}$ of $\{f_n\}$ converging uniformly on
$\overline{D}(a;r_1,r_2)$ to a bicomplex (holomorphic) function
$f$. That is, for some positive integer $n_0$, we have
$$\|f_{n_k}(w)-f(w)\|<1, \ \ \ \forall k\geq n_0,\text{ and } w\in \overline{D}(a;r_1,r_2).$$
Thus, if $M=max_{z\in \overline{D}(a;r_1,r_2)}\|f(w)\|$, then
$\|f_{n_k}(w)\|\leq 1+M,\ \ \forall w\in \overline{D}(a;r_1,r_2)$
and this is a contradiction.

\medskip
\noindent The above discussion permits to establish the following results.

\begin{theorem}
A family $\boldsymbol F$ of bicomplex holomorphic functions is
normal on the arbitrary domain $D$ if and only if $\boldsymbol F_{ei}$=
$\mathcal{P}_{i}(\boldsymbol F)$ is normal on $\mathcal{P}_{i}(D)$ for $i=1,2.$
\label{maintheo}
\end{theorem}

\begin{corollary}
If a family $\boldsymbol F$ of bicomplex holomorphic functions
is normal on a arbitrary domain $D\neq \mathcal{P}_1(D)\times_{e}\mathcal{P}_2(D)$, then $\boldsymbol F$ is
normal on the larger domain
$\mathcal{P}_1(D)\times_{e}\mathcal{P}_2(D)$.
\end{corollary}

\begin{corollary}
A family $\boldsymbol F$ of bicomplex holomorphic  functions is normal on a arbitrary domain $D$ if and only if $\boldsymbol F$ is normal at
each point of $D$.
\end{corollary}

\section{Bicomplex Montel Theorem from Montel Theorem of $\mathbb{C}^{2}$}

In this section, we want to show that it is possible to see the Bicomplex Montel Theorem (Theorem \ref{Montel}) as a particular case of the following Montel theorem of several complex variables (see \cite{26}).

\begin{theorem}
Let $D\subset\mathbb{C}^{n}$ be an open set and $\boldsymbol{F}\subset\mathcal{O}(D,\mathbb{C}^{n})$ be a family of holomorphic mappings.
Then the following are equivalent:
\begin{itemize}
\item[1.] The family $\boldsymbol{F}$ is locally uniformly bounded.
\item[2.] The family $\boldsymbol{F}$ is relatively compact in $\mathcal{O}(D,\mathbb{C}^{n})$.
\end{itemize}
\label{BMC2}
\end{theorem}

Since, $TH(D)\subset\mathcal{O}(D,\mathbb{C}^{2})$, we obtain directly the desired result using the fact that a family
$\boldsymbol{F}$ is relatively compact in $\mathcal{O}(D,\mathbb{C}^{n})$ a family $\boldsymbol{F}$ is relatively compact in $\mathcal{O}(D,\mathbb{C}^{n})$ if and only if $\boldsymbol{F}$ is a normal family. Recall that a family $\boldsymbol{F}$ is said
to be relatively compact if the family $\overline{\boldsymbol{F}}$ is compact (see \cite{33}).
Moreover, Theorem \ref{BMC2} will be proved for the specific
class $TH(D)$ instead of $\mathcal{O}(D,\mathbb{C}^{2})$ if we can show that $TH(D)$ is closed in
$\mathcal{O}(D,\mathbb{C}^{2})$ with the compact convergence topology. This is a direct consequence
of the following \textbf{Bicomplex Weierstrass Theorem}.

\begin{lemma}
Let $\{f_n\}$ be a sequence of bicomplex holomorphic functions which converges locally uniformly to a function $f$ on a
$\mathbb{T}$-disc $D(a_1+a_2\bold{i_2};r,r)$.
Then $f$ is bicomplex holomorphic in $D(a_1+a_2\bold{i_2};r,r)$.
\label{weierslem}
\end{lemma}
\emph{Proof}
Since $f_n(z_1+z_2\bold{i_2})$ is $\mathbb{T}$-holomorphic on $D(a_1+a_2\bold{i_2};r,r)$ $\forall n\in\mathbb{N}$, we have from Theorem \ref{theo4} that
$$(f_{ei})_n:\mathcal{P}_i(D(a_1+a_2\bold{i_2},r))\longrightarrow\mathbb{C}(\bold{i_1})$$ is holomorphic for $i=1,2$, $\forall n\in\mathbb{N}$.
Since $D(a_1+a_2\bold{i_2};r,r)$ is a bicomplex cartesian product, by the Weierstrass theorem of one complex variable, the sequence $(f_{ei})_n$ must converges locally uniformly
to the holomorphic function $f_{ei}$ on $D(\mathcal{P}_i(a_1+a_2\bold{i_2}),r)$ for $i=1,2$. Therefore, from Theorem \ref{theo5}, the function
$f(z_1+z_2\bold{i_2})=f_{e1}(z_1-z_2\bold{i_1})\bold{e_1}+f_{e2}(z_1+z_2\bold{i_1})\bold{e_2}$ is $\mathbb{T}$-holomorphic
on $D(a_1+a_2\bold{i_2};r,r)$.$\Box$

\begin{theorem}
\textbf{(Weierstrass)} Let $\{f_n\}$ be a sequence of bicomplex holomorphic functions on a domain $D$ which converges uniformly on compact subsets of $D$
to a function $f$. Then $f$ is bicomplex holomorphic in $D$.
\end{theorem}
\emph{Proof}
For an arbitrary $w_0\in D$, choose a $\mathbb{T}$-disc $D(w_0;r,r)\subset D$. Since $f_n(w)\rightarrow f(w)$ locally uniformly on $D$,
by Lemma \ref{weierslem}, $f$ is $\mathbb{T}$-holomorphic on $D$. As $w_0$ was arbitrary, $f(w)$ is $\mathbb{T}$-holomorphic on $D$.$\Box$

\section{A More General Definition of Normality}

To carry further the study of normal families of  bicomplex holomorphic functions particularly to consider the dynamics of
bicomplex holomorphic functions, we propose the following more general definition of normality.

\begin{definition}  The family $\boldsymbol F$ of bicomplex holomorphic functions defined on a domain $D \subseteq \mathbb{T}$ is said to be
\textbf{normal} in $D$ if every sequence in  $\boldsymbol F$
contains a subsequence which on compact subsets of $D$ either
converges uniformly to a limit function or converges uniformly to
$\infty$. $\boldsymbol F$ is said to be \textbf{normal at a point}
$z\in D$ if it is normal in some neighborhood of $z$ in $D$.
\label{compnormal}
\end{definition}
\begin{remark}
We say that the sequence $\{w_n\}$ of bicomplex numbers converges to
$\infty$ if and only if the norm $\{\|w_n\|\}$ congerges to
$\infty$.
\end{remark}

We note that our proofs of the Bicomplex Montel Theorem work in this situation too. However,
as for one complex variable, the converse of  Theorem \ref{Montel} will not remain valid with this more complete definition
of  the normality (see \cite{27}).

\begin{remark} Both situations in the last definition may occur simultaneously. For example, consider the family
$\{R^{\circ n}(w) \mid \ R(w)=w^2\mbox{ and  } n\in \mathbb N \}$ of bicomplex holomorphic functions on $\mathbb {T}$.
Then, by using the idempotent representation and results from one complex variable theory of normal families, we find
that this family is normal on $A\cup B$, where
$$A=\{w=w_1\eo +w_2\et:\ |w_1|<1, \ |w_2|<1\ \}$$
and
$$B=\{w=w_1\eo+w_2\et:\ |w_1|>1, \ |w_2|>1\ \}$$
On the set $A$, normality is under the first situation whereas on the set $B$ the normality is under the second situation.
\end{remark}

\begin{example}  Consider the family
$$\boldsymbol F =\{ f_{n}(w)=nw: w=z_1+z_2\bold{i_2}, \ \ n\in \mathbb Z \}.$$
Then $f_{n}(0)\rightarrow 0$, but $f_{n}(w)\rightarrow \infty$ for $w\neq 0$. It follows that $\boldsymbol F$ cannot be normal
in any domain containing the origin.
\end{example}

Now, let us prove that Theorem \ref{maintheo} is only true in one direction with this more general definition
of normality.

\begin{theorem}
Let $\boldsymbol F$ be a family of bicomplex holomorphic functions defined on a domain $D$.
If $\boldsymbol F_{ei}=\mathcal{P}_{i}(\boldsymbol F)$ is normal on $\mathcal{P}_{i}(D)$ for $i=1,2$ then $\boldsymbol F$
is normal on $D$.
\label{normal2}
\end{theorem}
\emph{Proof}
Suppose that $\boldsymbol F_{ei}=\mathcal{P}_{i}(\boldsymbol F)$ is normal on $\mathcal{P}_{i}(D)=D_i, \ i=1,2.$
We want to show that $\boldsymbol F$ is normal in $D$. Let $\{F_n\}$ be any sequence in $\boldsymbol F$ and
$K$ be any compact subset of $D$. Then $\{\mathcal{P}_1(F_n)\}=\{(f_n)_1\}$ is a sequence in
$\boldsymbol F_{e1}=\mathcal{P}_{1}(\boldsymbol F).$ Since $\boldsymbol F_{e1}=\mathcal{P}_{1}(\boldsymbol F)$ is normal in
$\mathcal{P}_1(D)$ then $\{(f_n)_1\}$ has a subsequence $\{(f_{n_k})_1 \}$ which converges
uniformly on $\mathcal{P}_1(K)$ to either a $\mathbb{C}(\bo)$-function or to $\infty$. Now, consider $\{F_{n_k}\}$ in $\boldsymbol F$.
Then $\{\mathcal{P}_2(F_{n_k})\}=\{(f_{n_k})_2\}$ is a sequence in
$\boldsymbol F_{e2}=\mathcal{P}_{2}(\boldsymbol F)$. Since $\boldsymbol F_{e2}=\mathcal{P}_{2}(\boldsymbol F)$ is normal in
$\mathcal{P}_1(D)$ then $\{(f_{n_k})_2\}$ has a subsequence $\{(f_{n_{k_l}})_2 \}$ which converges
uniformly on $\mathcal{P}_2(K)$ to either a $\mathbb{C}(\bo)$-function or to $\infty$. This implies that $\{(f_{n_{k_l}})_1\eo+(f_{n_{k_l}})_2\et\}$
is a subsequence of $\{\boldsymbol F_n\}$ which converges uniformly on $\mathcal{P}_1(K)\times_{e}\mathcal{P}_2(K)\supseteq K$ to a
bicomplex function or to $\infty$ showing that $\boldsymbol F$ is normal in $D$. $\Box$

\bigskip
\noindent Here is a counterexample for the other side.
\begin{example}
Let $X_1$ and $X_2$ be domains in $\mathbb{C}(\bold{i_1})$ containing the origin. Let $D = (X_1\times_e X_2) - \{0\}$.
Then $D$ is not a bicomplex cartesian domain because $D\neq \mathcal{P}_1(D)\times_e \mathcal{P}_2(D)$. Now the family

$$\boldsymbol F =\{ nw: w=z_1+z_2\bold{i_2}, \ \ n\in \mathbb N \}$$
is normal in the domain $D$ (by the proposed definition of normality as above) but $\boldsymbol F_{ei}=\mathcal{P}_{i}(\boldsymbol F)$ is not normal
in $\mathcal{P}_i(D), \ i=1,2$ as it contains the origin.
\end{example}

Moreover, the next examples show that the converse of Theorem \ref{normal2} is not true even if the domain
$D$ is a bicomplex cartesian product.

\begin{example}
Consider the family
$$\mathcal F=\{nz: z\in \mathbb{C}(\bo), \ \ n\in \mathbb Z \}$$
on $\mathbb{C}(\bo)$.
Then $\mathcal F$ is normal on the punctured disc $D(0,1)-\{0\}\subset \mathbb{C}(\bo)$ but not normal on the disc $D(0,1)\subset \mathbb{C}(\bo)$.
However, the bicomplex family $$\boldsymbol F:=\boldsymbol F_{e1}\eo+\boldsymbol F_{e2}\et=\{ nw: w=z_1+z_2\bold{i_2}, \ \ n\in \mathbb N \}$$
where $\boldsymbol F_{e1}=\mathcal F$ is normal in the following bicomplex cartesian product:
$$(D(0,1)-\{0\})\times_e D(0,1)$$
since the limit function is identically infinite.
\end{example}

\begin{example}
Consider the family
$$\mathcal F=\{R^{\circ n}(z) \mid \ R(z)=z^2 \mbox{ and  } n\in \mathbb N\}$$
on $\mathbb{C}(\bo)$.
Then $\mathcal F$ is normal on $D_1=\{z:|z|>1\}\subset \mathbb{C}(\bo)$ where here the limit function is identically infinite,
but not normal on $\mathbb{C}(\bo)$ since $\{|z|=1\}\subset \mathbb{C}(\bo)$.
However, the bicomplex family $$\boldsymbol F:=\boldsymbol F_{e1}\eo+\boldsymbol F_{e2}\et=\{R^{\circ n}(w) \mid \ R(w)=w^2 \mbox{ and  } n\in \mathbb N \}$$
where $\boldsymbol F_{e1}=\boldsymbol F_{e2}=\mathcal F$, is normal in the following bicomplex cartesian product:
$$D_1\times_e \mathbb{C}(\bo)$$
since the limit function is identically infinite.
\label{maincounter}
\end{example}

\section{Foundation of Bicomplex Dynamics: Fatou and Julia Sets for Polynomials}

We conclude this article with the following general definition of \textbf{Fatou} and \textbf{Julia} sets for bicomplex polynomials.
\begin{definition}
Let $P(\zeta)$ be a bicomplex polynomial. We define the bicomplex Julia set for $P$ as
\begin{equation}
\mathcal{J}_{2}(P)=\{\zeta\in\mathbb{T} \mid \{P^{\circ n}(\zeta)\}\mbox{ is not normal}\}
\label{fund}
\end{equation}
and the bicomplex Fatou (or stable) set as
\begin{equation}
\mathcal{F}_{2}(P)=\mathbb{T}-\mathcal{J}_{2}(P).
\end{equation}
\end{definition}

\noindent Hence, for each point $\zeta\in \mathcal{F}_{2}(P)$, there is a
neighborhood $N_{\zeta}$ in which $\{P^{\circ n}(\zeta)\}$ is a
normal family. Therefore, $\mathcal{F}_{2}(P)$ is an open set, the connected
components of which are the maximal domains of normality of
$\{P^{\circ n}(\zeta)\}$, and $\mathcal{J}_{2}(P)$ is a closed set.

From Theorem \ref{normal2}, we obtain the following inclusion:

$\mathcal{J}_{2}(P) \subset \mbox{\large \{ }z_1+z_2\bt\in\mathbb{T}\mid \{[\mathcal{P}_1(P)]^{\circ n}(z_1-z_2\bo)\}\}\mbox{ or}$
\begin{eqnarray}
        & & \{[\mathcal{P}_2(P)]^{\circ n}(z_1+z_2\bo)\}\mbox{ is not normal \large \}}\label{inc}\\
        &=& [\mathcal{J}_1(\mathcal{P}_1(P))\times_e \mathbb{C}(\bo)]\cup [\mathbb{C}(\bo)\times_e \mathcal{J}_1(\mathcal{P}_2(P))].
\end{eqnarray}

However, from Example \ref{maincounter}, we know that (\ref{inc})
cannot be transformed into an equality. In fact, to obtain a
characterization of bicomplex Julia sets in terms of one variable
dynamics we need to use the concept of filled-in Julia set. As for
the complex case, the bicomplex \textbf{filled-in Julia} set
$K_2(P)$ of a polynomial $P$ is defined as the set of all points
$\zeta$ of dynamical space that have bounded orbit with respect to
$P$, that is to say
\begin{equation}
\mathcal{K}_2(P)=\{\zeta\in\mathbb{T} \mid \{P^{\circ n}(\zeta)\}\nrightarrow\infty\mbox{ as }n\rightarrow\infty \}.
\end{equation}
\noindent We remark that $\mathcal{K}_2(P)$ is a closed set.
\newpage

As for the classical case (see \cite{29}, P.65), we need to
consider polynomials of degree $d\geq 2$ to be able to see a
bicomplex Julia set as the boundary of a bicomplex filled-in Julia
set. In fact, to decompose $P(w)$ in terms of two complex
polynomials of degree $d\geq 2$, we must also consider non-degenerate
bicomplex polynomials of the form
$P(w)=a_dw^d+a_{d-1}w^{d-1}+...+a_0$ where $a_d\notin
\mathcal{NC}$.

\smallskip
\noindent Under these specifications, we have the following
result.

\begin{theorem}
Let $P(\zeta)$ be a non-degenerate bicomplex polynomials of degree $d\geq 2$. Then,
\begin{equation}
\partial \mathcal{K}_2(P)=\mathcal{J}_{2}(P).
\end{equation}
\label{filled}
\end{theorem}

Now, using the concept of normality in terms of \textbf{Definition \ref{basicnormal}}, we obtain the following characterization of
$\mathcal{K}_{2}(P)-\mathcal{J}_{2}(P)$:
\begin{equation}
\mathcal{K}_{2}(P)-\mathcal{J}_{2}(P)=\{\zeta\in\mathbb{T} \mid \{P^{\circ n}(\zeta)\}\mbox{ is \textit{normal}}\}.
\end{equation}

Moreover, using the idempotent representation, it is easy to see
that the bicomplex filled-in Julia set $\mathcal{K}_2(P)$ can be
expressed in terms of the two filled-in Julia sets in the plane. More
specifically,
\begin{equation}
\mathcal{K}_2(P)=\mathcal{K}_1(\mathcal{P}_1(P))\times_e \mathcal{K}_1(\mathcal{P}_2(P)).
\end{equation}

Hence, since
$\partial[\mathcal{K}_1(\mathcal{P}_1(P))\times_e \mathcal{K}_1(\mathcal{P}_2(P))]=[\partial \mathcal{K}_1(\mathcal{P}_1(P))\times_e \mathcal{K}_1(\mathcal{P}_2(P))]\cup [\mathcal{K}_1(\mathcal{P}_1(P))\times_e \partial \mathcal{K}_1(\mathcal{P}_2(P))]$,
we have the following characterization of the bicomplex Julia set $\mathcal{J}_{2}(P)$ in terms of one complex variable dynamics.
\begin{theorem}
Let $P(\zeta)$ be a non-degenerate bicomplex polynomials of degree $d\geq 2$. Then,
\begin{equation}
\mathcal{J}_{2}(P)=[\mathcal{J}_1(\mathcal{P}_1(P))\times_e \mathcal{K}_1(\mathcal{P}_2(P))]\cup [\mathcal{K}_1(\mathcal{P}_1(P))\times_e \mathcal{J}_1(\mathcal{P}_2(P))] \label{put}.
\end{equation}
\label{grabache}
\end{theorem}

\begin{example}
Consider the bicomplex polynomial: $$P(w)=w^2.$$ We can verify that $\mathcal{P}_k(w^2)=z^2$ for $k=1,2$.
In the complex plane (in $\bo$), it is well known that $\mathcal{K}_1(z^2)=\{z:|z|<1\}$ and
$\mathcal{J}_1(z^2)=\{z:|z|=1\}$ (see \cite{30} or \cite{27}).
Hence, using Theorem \ref{grabache}, we obtain that
\begin{equation}
\mathcal{J}_{2}(P(w))=[\{z:|z|=1\}\times_e \{z:|z|<1\}]\cup [\{z:|z|<1\}\times_e \{z:|z|=1\}].
\end{equation}
\end{example}

\begin{remark}
In the particular case of the bicomplex quadratic polynomial
\begin{equation}
P_c(\zeta)={\zeta}^2+c,
\end{equation}
the \textbf{fundamental} definition of the bicomplex Julia set of this
article (see \ref{fund}) coincides with the definition, using boundary of bicomplex filled-in
Julia set, introduced by D. Rochon in \cite{16,28} (see Theorem \ref{filled}). Moreover, using some distance estimation formulas
that can be used to ray traced slices of bicomplex filled-in Julia
sets in dimension three (see \cite{15}), we obtain some
visual examples (see Fig. 1, 2, 3 and 4) of bicomplex Julia sets
$\mathcal{K}_2(P_c)$ for the specific slice $\bj=0$.
\end{remark}

\begin{example}
Consider the bicomplex polynomial: $$P(w)=w^2+0.27.$$ We can verify that $\mathcal{P}_k(w^2+0.27)=z^2+0.27$ for $k=1,2$.
In the complex plane (in $\bo$), it is well known that $\mathcal{K}_1(z^2+0.27)=\mathcal{J}_1(z^2+0.27)$ is a Cantor set.
We shall denote such Cantor set by $\mathcal{C}_{0.27}$.
Hence, using Theorem \ref{grabache}, we obtain that
\begin{equation}
\mathcal{J}_{2}(P)=\mathcal{C}_{0.27}\times_e \mathcal{C}_{0.27} \ \ \mbox{(Fig. 1)}.
\end{equation}
\end{example}

\begin{example}
Consider the bicomplex polynomial: $$P(w)=w^2-1.754878.$$ We can verify that $\mathcal{P}_k(w^2-1.754878)=z^2-1.754878$ for $k=1,2$.
In the complex plane (in $\bo$), it is well known that $\mathcal{K}_1(z^2-1.754878)$ is the so-called \textbf{Airplane} (see \cite{29}, P.129).
We shall denote this set by $\mathcal{A}$.
Hence, using Theorem \ref{grabache}, we obtain that
\begin{equation}
\mathcal{J}_{2}(P)=[\partial \mathcal{A}\times_e \mathcal{A}]\cup [\mathcal{A}\times_e \partial \mathcal{A}] \ \ \mbox{(Fig. 2)}.
\end{equation}
\end{example}

\begin{figure}
\centering
%c=(0.27)e_1+(0.27)e_2
\includegraphics[width=3.9cm]{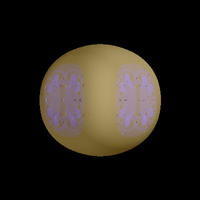}
\includegraphics[width=3.9cm]{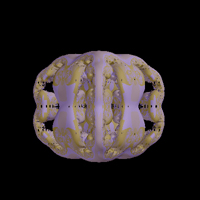}
\includegraphics[width=3.9cm]{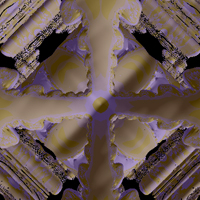}
\caption{Three Explorations of $\mathcal{K}_2(P_c)$ for $c=(0.27)\eo+(0.27)\et$}
\end{figure}

\begin{figure}
\centering
%c=(-1.754878)e_1+(-1.754878)e_2
\includegraphics[width=3.9cm]{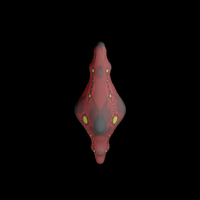}
\includegraphics[width=3.9cm]{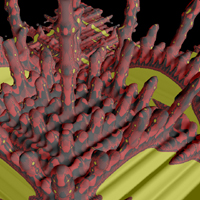}
\includegraphics[width=3.9cm]{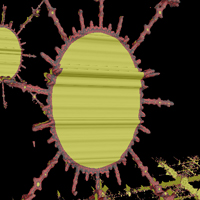}
\caption{Three Explorations of $\mathcal{K}_2(P_c)$ for $c=(-1.754878)\eo+(-1.754878)\et$}
\end{figure}

\newpage
\begin{example}
Consider the bicomplex polynomial: $$P(w)=w^2+[(0.26)\eo+(-1.754878)\et].$$ We can verify that
$\mathcal{P}_1(w^2+[(0.26)\eo+(-1.754878)\et])=z^2+0.26$ and
$\mathcal{P}_2(w^2+[(0.26)\eo+(-1.754878)\et])=z^2-1.754878$.
Hence, using Theorem \ref{grabache}, we obtain that
\begin{equation}
\mathcal{J}_{2}(P)=\mathcal{C}_{0.26}\times_e \mathcal{A} \ \ \mbox{(Fig. 3)}.
\end{equation}
\end{example}

\begin{figure}
\centering
%c=(0.26)e_1+(-1.754878)e_2
\includegraphics[width=3.9cm]{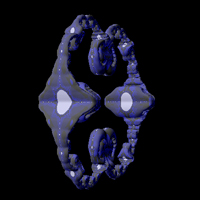}
\includegraphics[width=3.9cm]{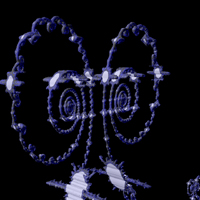}
\includegraphics[width=3.9cm]{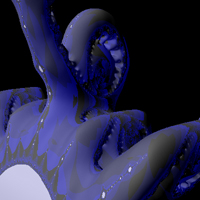}
\caption{Three Explorations of $\mathcal{K}_2(P_c)$ for $c=(0.26)\eo+(-1.754878)\et$}
\end{figure}

\begin{example}
Consider the bicomplex polynomial: $$P(w)=w^2+[(-0.123+0.745\bo)\eo+(-0.391-0.587\bo)\et].$$ We can verify that $\mathcal{P}_1(w^2+[(-0.123+0.745\bo)\eo+(-0.391-0.587\bo)\et])=z^2+(-0.123+0.745\bo)$ and
$\mathcal{P}_2(w^2+[(-0.123+0.745\bo)\eo+(-0.391-0.587\bo)\et])=z^2+(-0.391-0.587\bo)$.
In the complex plane (in $\bo$), it is well known that $\mathcal{D}:=\mathcal{K}_1(z^2+(-0.123+0.745\bo))$ is the so-called \textbf{Douady's Rabbit}
and $\mathcal{S}:=\mathcal{K}_1(z^2+(-0.391-0.587\bo))$ is a \textbf{Siegel Disk}.
Hence, using Theorem \ref{grabache}, we obtain that
\begin{equation}
\mathcal{J}_{2}(P)=[\partial \mathcal{D}\times_e \mathcal{S}]\cup [\mathcal{D}\times_e \partial \mathcal{S}] \ \ \mbox{(Fig. 4)}.
\end{equation}
\end{example}

\begin{figure}
\centering
%c=(-0.123+0.745i_1)e_1+(-0.391-0.587i_1)e_2
\includegraphics[width=6cm]{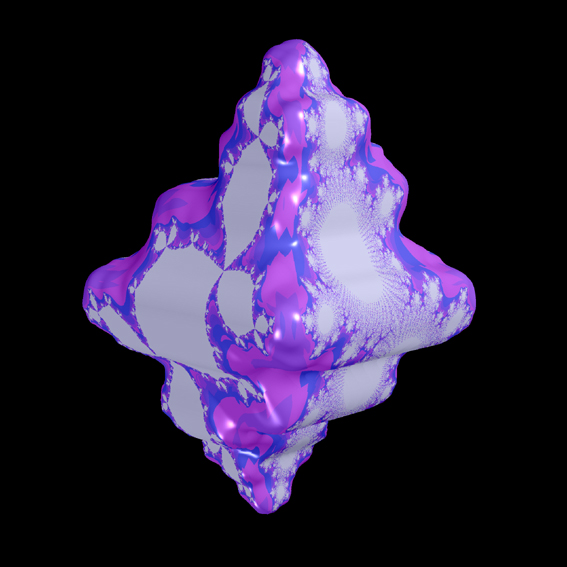}
\includegraphics[width=6cm]{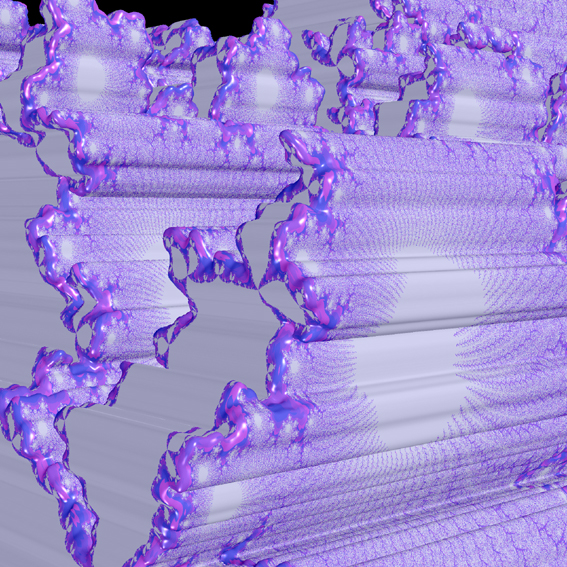}
\caption{Two Explorations of $\mathcal{K}_2(P_c)$ for $$c=(-0.123+0.745\bo)\eo+(-0.391-0.587\bo)\et$$}
\end{figure}

\begin{remark} By using the definition of the bicomplex Fatou set as the complement of the bicomplex Julia set (\ref{put}) leads us to
characterize the bicomplex Fatou set of non-degenerate bicomplex polynomials of degree $d\geq 2$ as
$$\mathcal{F}_2(P)=[\mathcal{F}_1(\mathcal{P}_1(P))\times_e \mathcal{F}_1(\mathcal{P}_2(P))]\cup [\mathcal{F}_1(\mathcal{P}_1(P))_{\infty}\times_e \mathcal{J}_1(\mathcal{P}_2(P))]$$
\begin{equation}
\cup [\mathcal{J}_1(\mathcal{P}_1(P))\times_e \mathcal{F}_1(\mathcal{P}_2(P))_{\infty}]\label{do}
\end{equation}
where $\mathcal{F}_1(\mathcal{P}_i(P))_{\infty}, \ i=1,2$ denotes the unbounded component of the Fatou set of the projections of $P.$
Moreover, from (\ref{do}) it follows that if a family $\boldsymbol P$ of bicomplex polynomials is normal in an unbounded domain $D$ of $\mathbb{T}$ then
at least one of the projections $\mathcal{P}_{i}(\boldsymbol P)$ is normal on the corresponding unbounded domain $\mathcal{P}_{i}(D)$, $i=1,2.$
\end{remark}

\begin{example}
Consider the bicomplex polynomial: $$P(w)=w^2.$$
In the complex plane (in $\bo$), it is well known that $$\mathcal{F}_1(z^2)=\{z:|z|<1\}\cup\{z:|z|>1\}$$ and
$\mathcal{F}_1(z^2)_{\infty}=\{z:|z|>1\}$.
Hence, using (\ref{do}), we obtain,

$$\mathcal{F}_2(P)=[(\{z:|z|<1\}\cup\{z:|z|>1\})\times_e (\{z:|z|<1\}\cup\{z:|z|>1\})]$$
$$\cup[\{z:|z|>1\}\times_e \{z:|z|=1\}]$$
\begin{equation}
\cup[\{z:|z|=1\}\times_e \{z:|z|>1\}].
\end{equation}
\end{example}

As a direct consequence of Theorem \ref{grabache}, the bicomplex Julia set $\mathcal{J}_{2}(P)$ is completely invariant under the substitution
$(w,P(w))$ when $P$ is a non-degenerate bicomplex polynomial of degree $d\geq 2$. The next theorem proves this
result in general.

\begin{definition}
Let $f(z_1+z_2\bold{i_2})=f_{e1}(z_1-z_2\bold{i_1})\eo+f_{e1}(z_1+z_2\bold{i_1})\et:D\longrightarrow\mathbb{T}$ be a bicomplex function. The function $f$ is said to be strongly non-constant on $D$ if $f_{ei}$ is non-constant on $\mathcal{P}_{i}(D)$ for $i=1,2$.
\end{definition}

\begin{theorem}
Let $f$ be an entire strongly non-constant $\mathbb{T}$-holomorphic function and
$$\mathcal{J}_2(f):=\{\zeta\in\mathbb{T} \mid \{f^{\circ n}(\zeta)\}\mbox{ is not normal}\}.$$
Then,\\
1. If a point $w_0\in\mathcal{J}_2(f)$, then $f(w_0)\in\mathcal{J}_2(f)$;\\
2. If $w_0\in\mathcal{J}_2(f)$ and $w_1$ is a point such that $f(w_1)=w_0$ then $w_1\in\mathcal{J}_2(f)$.
\end{theorem}
\emph{Proof}
Since our definition of normality (Def. \ref{compnormal}) is analogous to the related notion in the complex plane,
the proof of the theorem is same as the proof of the corresponding result in the plane (see \cite{31}, Theorem 2.17 and 2.18). \textbf{Note:} The condition for $f$ to be strongly non-constant is needed in the proof
of (2.) to be able to use the open mapping theorem in each components.$\Box$

\section*{Acknowledgments}
The research of D.R. is partly supported by grants from CRSNG of Canada and FQRNT of Qu{\'e}bec.

\end{document}